\documentclass[12pt,leqno]{article}
\usepackage{amsmath,amssymb,amscd,latexsym}
\newcommand{\C}{{\mathbb{C}}}
\newcommand{\R}{{\mathbb{R}}}
\newcommand{\ad}{\mathrm{ad}}
\newcommand{\Ad}{\mathrm{Ad}}
\newcommand{\even}{\mathrm{even}}
\newcommand{\fin}{\mathrm{fin}}
\newcommand{\id}{\mathrm{id}}
\newcommand{\imm}{\mathrm{im}\,}
\renewcommand{\mod}{\;\mathrm{mod}\;}
\newcommand{\sign}{\mathrm{sign}\,}
\newcommand{\esl}{\mathfrak{sl}}
\newcommand{\so}{\mathfrak{so}}
\newcommand{\vol}{\mathrm{vol}}
\newcommand{\End}{\mathrm{End}\,}
\newcommand{\Hom}{\mathrm{Hom}}
\newcommand{\Ker}{\mathrm{Ker}\,}
\newcommand{\Ah}{{\mathcal A}}
\newcommand{\Ch}{{\mathcal C}}

\newcommand{\Fh}{{\mathcal F}}
\newcommand{\Hh}{{\mathcal H}}
\newcommand{\Kh}{{\mathcal K}}
\newcommand{\ea}{\mathfrak{a}}
\newcommand{\eg}{{\mathfrak{g}}}
\newcommand{\ek}{\mathfrak{k}}

\newcommand{\ep}{\mathfrak{p}}
\newcommand{\eU}{\mathfrak{U}}
\newcommand{\obeta}{\overline{\beta}}
\newcommand{\ove}{\overline{e}}
\newcommand{\oeta}{\overline{\eta}}

\newcommand{\oBox}{\overline{\square}}
\newcommand{\oH}{\bar{H}}
\newcommand{\opartial}{\bar{\partial}}
\newcommand{\ou}{\overline{u}}
\newcommand{\ohne}{\setminus}
\newcommand{\silo}{\stackrel{\sim}{\longrightarrow}}
\newcommand{\tei}{\, | \,}

\newcommand{\hullet}{\raisebox{0.05cm}{$\scriptscriptstyle \bullet$}}
\newcommand{\verk}{\mbox{\scriptsize $\,\circ\,$}}
\newcommand{\halb}{\frac{1}{2}}

\newtheorem{theorem}{Theorem}[section]
\newtheorem{prop}[theorem]{Proposition}
\newtheorem{lemma}{Lemma}[theorem]
\newenvironment{defn}{\noindent \textbf{Definition}}{}
\newenvironment{remark}{\noindent \textbf{Remark}}{}
\newtheorem{remarks}[theorem]{Remark}
\newtheorem{punkt}[theorem]{$\!\!$}
\newenvironment{proof}{\noindent {\bf Proof}}{\mbox{}\hfill$\Box$}
\begin{document}

%%%%%%%%%%%%%%%%%%%%%%%%%%%%%%%%%%%%%%%%%%%%%%%%%%%%%%%%%%%%%%%%

\title{Real polarizable Hodge structures arising from foliations}

\author{by Christopher Deninger and Wilhelm Singhof}
\date{\ }
%\address{Mathematisches Institut, WWU M\"unster, Einsteinstr. 62, 48149 M\"un\-ster, Germany}
%\curraddr{...}
%\email{deninge@math.uni-muenster.de}
%\thanks{...}

%    General info
%\subjclass{37Cxx, 37C27, 53C12}
%\date{January 1, 1994 and, in revised form, June 22, 1994.}

%\dedicatory{...}

%\keywords{Flow, foliation, trace formula, leafwise cohomology}

%\begin{abstract}
%We prove a dynamical Lefschetz trace formula for certain flows respecting a one-codimensional foliation.
%\end{abstract}
\maketitle
%\hspace*{\fill} {\footnotesize preliminary uncorrected version}

%%%%%%%%%%%%%%%%%%%%%%%%%%%%%%%%%%%%%%%%%%%%%%%%%%%%%%%%%%%%%%%%
%\input{intro}

\section*{Introduction}
\label{sec:0}

Consider a smooth foliation $\Fh$ on a connected closed manifold $X$. The leafwise cohomology of $\Fh$ is the cohomology of the de Rham complex of smooth forms along the leaves. Even its maximal Hausdorff quotient, the reduced leafwise cohomology $\oH^{\hullet}_{\Fh} (X)$, is a quite subtle invariant of the foliation. In particular it can be infinite dimensional even in degree one.

For general foliations there does not seem to be a good structure theory for $\oH^{\hullet}_{\Fh} (X)$. For Riemannian foliations, the situation is different however. This is mainly due to the smooth Hodge decomposition theorem of \'Alvarez L\'opez and Kordyukov \cite{AK}. This theorem expresses $\oH^{\hullet}_{\Fh} (X)$ as the space of smooth global forms on $X$ along $\Fh$ whose restrictions to the leaves are harmonic. 

Given such a Hodge theorem it is a very natural question whether Hodge--K\"ahler--Lefschetz theory can be generalized to Riemannian $\Fh$ if the leaves are K\"ahler manifolds. Of course the K\"ahler structures should vary smoothly in the transverse direction. In this paper we show that this program can be carried through. In fact it can be done quite easily and the first four sections essentially just formulate the resulting assertions, the proofs being mostly consequences of the classical K\"ahler--Lefschetz identities on non-compact K\"ahler manifolds. The results are the same as the usual ones except that our Hodge structures can be infinite dimensional.

As in the classical case the existence of a K\"ahler structure on a Riemannian foliation imposes a number of restrictions on the structure of $\oH^{\hullet}_{\Fh} (X)$, in particular those implied by the hard Lefschetz theorem.

Incidentally the classical case is contained in the foliation formalism by taking for $\Fh$ the foliation which consists of the leaf $X$ only. In this case $\oH^{\hullet}_{\Fh} (X)$ becomes the ordinary cohomology of $X$ with real coefficients.

As an application of the K\"ahler formalism we transfer to $\oH^{\hullet}_{\Fh} (X)$ a result of Serre which establishes an analogue of the Weil conjectures for certain induced ``Frobenius''-maps on cohomology, see \ref{t411}.

In \cite{LL} Looijenga and Lunts attach a Lie algebra to each compact K\"ahler manifold using the Hodge--K\"ahler--Lefschetz formalism. We extend their construction to K\"ahler--Riemann foliations. In this way one obtains a possibly infinite dimensional Lie algebra attached to a K\"ahler--Riemann foliation. In proposition \ref{t52} we determine its structure in the case where $\Fh$ is a foliation by complex surfaces with a dense leaf. 

Finally in section \ref{sec:6} we discuss a class of examples suggested to us by E. Ghys. They are the foliations obtained by projecting the foliations $S_1 \times \{ pt \}$ and $\{ pt \} \times S_2$ to the quotient $X = \Gamma \ohne (S_1 \times S_2)$. Here $S_i = G_i / K_i$ is a bounded symmetric domain for $i = 1,2$ and the cocompact lattice $\Gamma \subset G_1 \times G_2$ is supposed to act without fixed points on $S_1 \times S_2$. The case of an irreducible lattice is particularly interesting because then the leaves are dense.

For these and more general K\"ahler--Riemann foliations we express the reduced leafwise cohomology in representation theoretic terms using $(\eg , K)$-cohomology. The final result gives a description of the Hodge $(p,q)$-decompo\-sition of $\oH^{\hullet}_{\Fh} (X)$ in these terms. The cohomology $\oH^{\hullet}_{\Fh} (X)$ in these examples is quite interesting and deserves further study.

In the article we have chosen to work with trivial coefficients only. More generally one may consider a Riemannian vector bundle $V$ with a Riemannian flat connection along $\Fh$. Using $V$-valued forms along $\Fh$ one defines a version of reduced leafwise cohomology with coefficients in $V$. All our results transfer to this setting. This is of interest for a possible generalization of the theory to variations of (polarizable) Hodge structures coming from foliated morphisms. 

Although it is not apparent, the original motivation for our paper arose from the number theoretical considerations in \cite{D} \S\,3 together with the later observation that the cohomology in question should be analogous to reduced leafwise cohomology \cite{D2} \S\,4.

We would like to thank Don Blasius for discussions about $(\eg , K)$-cohomology and E. Ghys for a letter pointing out the examples in section \ref{sec:6}.
\section{Complex structures on foliations}
\label{sec:1}
Let $\Fh$ be a smooth foliation of the manifold $X$. By $T\Fh$ we denote the subbundle of the tangent bundle $TX$ whose fibers $T_x \Fh$ are the tangent spaces of the leaves and by $T^* \Fh$ the dual bundle of $T\Fh$. Let
\[
\Ah^r_{\Fh} (X) := \Gamma (X , \Lambda^r T^* \Fh \otimes_{\R} \C)
\]
be the space of complex valued $r$-forms along the leaves.

If the dimension of the leaves is even, we define an \textit{almost complex structure} on $\Fh$ to be a smooth complex structure $J$ on $T\Fh$. That is, $J$ is a smooth real vector bundle automorphism of $T\Fh$ with $J^2 = - \id$. Given such an almost complex structure, $T\Fh$ becomes a complex vector bundle which we denote by $T^c \Fh$. The bundle $T^* \Fh \otimes_{\R} \C$, whose fiber over $x$ is the space of $\R$-linear maps $T_x \Fh \to \C$ then decomposes as
\[
  T^* \Fh \otimes_{\R} \C = T^{1,0} \Fh \oplus T^{0,1} \Fh
\]
where the fiber $T^{1,0}_x \Fh$ consists of the $\C$-linear maps $T^c_x \Fh \to \C$ and $T^{0,1}_x \Fh$ of the anti-linear ones. This leads to a decomposition of exterior algebras
\[
  \Lambda^{\hullet} T^* \Fh \otimes_{\R} \C = \Lambda^{\hullet} T^{1,0} \Fh \otimes_{\C} \Lambda^{\hullet} T^{0,1} \Fh \; .
\]
We obtain, with the notation 
\[
  \Ah^{p,q}_{\Fh} (X) := \Gamma ( X , \Lambda^p T^{1,0} \Fh \otimes \Lambda^q T^{0,1} \Fh) \; ,
\]
a bigrading
\[
  \Ah^{\hullet}_{\Fh} (X) = \bigoplus_{p,q} \Ah^{p,q}_{\Fh} (X) \; .
\]
An almost complex structure $J$ on $\Fh$ defines, for each leaf $\ell$, by restriction an almost complex structure $J_{\ell} : T\ell \to T\ell$ on the manifold $\ell$. We have correspondingly a bigrading
\[
  \Ah^{\hullet} (\ell) = \bigoplus_{p,q} \Ah^{p,q} (\ell)
\]
of complex valued differential forms on $\ell$. The following fact is obvious:

\begin{punkt}
  \label{t11} \rm
An element $\alpha \in \Ah^r_{\Fh} (X)$ belongs to $\Ah^{p,q}_{\Fh} (X)$ if and only if the restrictions $\alpha \tei \ell$ belong to $\Ah^{p,q} (\ell)$ for all leaves $\ell$.

If all the $J_{\ell}$ are integrable (which is the case if and only if all the leaves are complex manifolds) we say that $J$ is a \textit{complex structure} on $\Fh$. 

The exterior derivative along the leaves,
\[
  d_{\Fh} : \Ah^r_{\Fh} (X) \longrightarrow \Ah^{r+1}_{\Fh} (X) \; ,
\]
has the property that
\[
  (d_{\Fh} \alpha) \tei \ell = d_{\ell} (\alpha \tei l)
\]
where $d_{\ell}$ is the exterior derivative on the manifold $\ell$.
\end{punkt}

Hence we conclude from \ref{t11}:

\begin{punkt}
  \label{t12}
\rm If $J$ is a complex structure on $\Fh$, we obtain a decomposition
\[
  d_{\Fh} = \partial_{\Fh} + \opartial_{\Fh}
\]
where $\partial_{\Fh}$ has bidegree $(1,0)$ and $\opartial_{\Fh}$ has bidegree $(0,1)$. With the obvious notations, we have
\begin{eqnarray*}
  (\partial_{\Fh} \alpha) \tei \ell & = & \partial_{\ell} (\alpha \tei \ell) \; , \\
(\opartial_{\Fh} \alpha) \tei \ell & = & \opartial_{\ell} (\alpha \tei \ell) \; .
\end{eqnarray*}
\end{punkt}

\begin{punkt}
  \label{t13} \rm
Assume that $\dim X = a$. A $2g$-dimensional foliation $\Fh$ together with a complex structure can be described equivalently as follows. It is given by a maximal atlas of charts:
\[
  \phi : U \silo V \times W \; \mbox{with} \; U \subset X , V \subset \C^g , W \subset \R^{a-2g} \; \mbox{all open}
\]
such that the coordinate changes take the form:
\[
  \phi_2 \verk \phi^{-1}_1 (z,w) = (g_1 (z,w) , g_2 (w)) \; .
\]
Here $g_1$ and $g_2$ are $C^{\infty}$ and for every $w$ the function $z \mapsto g_1 (z,w)$ is holomorphic.
\end{punkt}

It is clear that such an atlas determines a foliation with a complex structure. The converse follows from the Newlander--Nirenberg theorem with parameters stated on p. 404 of \cite{NN}.

\begin{punkt}
  \label{t14} \rm
For a smooth foliation $\Fh$ let $\Ch_{\Fh}$ denote the sheaf of germs of smooth complex valued functions which are locally constant on the leaves. Let $\Ah^r_{\Fh}$ denote the sheaf of smooth sections of $\Lambda^r T^* \Fh \otimes_{\R} \C$. The primitive in a star shaped domain constructed in the standard proof of the Poincar\'e Lemma depends smoothly on parameters. There follows the well known fact that 
\[
  0 \longrightarrow \Ch_{\Fh} \longrightarrow \Ah^0_{\Fh} \overset{d_{\Fh}}{\longrightarrow} \Ah^1_{\Fh} \longrightarrow \ldots
\]
is a fine resolution of $\Ch_{\Fh}$. In particular we get an isomorphism
\[
  H^r (X , \Ch_{\Fh}) = H^r (\Ah^{\hullet}_{\Fh} (X) , d_{\Fh}) \; .
\]
The cohomology is called the leafwise cohomology of $\Fh$. Its maximal Hausdorff-quotient with respect to the natural Fr\'echet topologies on the $\Ah^r_{\Fh} (X)$'s is called the reduced leafwise cohomology of $\Fh$ with complex coefficients:
\[
  \oH^r_{\Fh} (X , \C) := \oH^r (X , \Ch_{\Fh}) = \oH^r (\Ah^{\hullet}_{\Fh} (X) , d_{\Fh}) \; .
\]
Using real valued forms, one defines in the same way the reduced leafwise cohomology $\oH^r_{\Fh} (X)$ of $\Fh$ with real coefficients.

If $\Fh$ carries a complex structure, let $\Omega^p_{\Fh}$ denote the sheaf of germs of smooth complex valued $p$-forms along the leaves which are holomorphic on the leaves. Let $\Ah^{pq}_{\Fh}$ denote the sheaf of smooth sections of $\Lambda^p T^{1,0} \Fh \otimes \Lambda^q T^{0,1} \Fh$. Using \ref{t13} a close inspection of the proof of the $\opartial$-Poincar\'e Lemma in \cite{GH} p. 25 shows that the primitive constructed there depends smoothly on parameters. If follows that
\[
0 \longrightarrow \Omega^p_{\Fh} \longrightarrow \Ah^{p,0}_{\Fh} \xrightarrow{\opartial_{\Fh}} \Ah^{p,1}_{\Fh} \longrightarrow \ldots
\]
is a fine resolution of $\Omega^p_{\Fh}$. Hence there are natural isomorphisms
\[
  H^q (X , \Omega^p_{\Fh}) \cong H^q (\Ah^{p,\cdot}_{\Fh} (X) , \opartial_{\Fh}) \; .
\]
The reduced versions of these leafwise Dolbeault cohomology groups are again defined by passing to the maximal Hausdorff-quotient:
\[
  \oH^q (X , \Omega^p_{\Fh}) := \oH^q (\Ah^{p,\cdot}_{\Fh} (X) , \opartial_{\Fh}) \; .
\]
\end{punkt}
\section{K\"ahler foliations}
\label{sec:2}

We consider a $2g$-dimensional foliation $\Fh$ on $X$ with a complex structure on $\Fh$. Let $h$ be an Hermitian metric on the bundle $T^c \Fh$. Generalizing the usual notation, we write
\[
h = S - 2i \omega_{\Fh}
\]
with real forms $S$ and $\omega_{\Fh}$. Then $S$ is a Riemannian metric on the bundle $T\Fh$, and $\omega_{\Fh} \in \Ah^{1,1}_{\Fh} (X)$. The leafwise $*$-operator
\[
*_{\Fh} :\Ah^{p,q}_{\Fh} (X) \longrightarrow \Ah^{g-q, g-p}_{\Fh} (X)
\]
is the pointwise operator defined by $(*_{\Fh} \alpha) \tei \ell = *_{\ell} (\alpha \tei \ell)$ where $*_{\ell}$ is the $*$-operator on the manifold $\ell$ endowed with the Riemannian metric obtained by restricting $S$ and with the orientation coming from the complex structure. We introduce the first order differential operators
\begin{eqnarray*}
  d^*_{\Fh} & := & -*_{\Fh} d_{\Fh} \, *_{\Fh} \; , \\
\partial^*_{\Fh} & := & - *_{\Fh} \opartial_{\Fh} \, *_{\Fh} \; , \\
\opartial^*_{\Fh} & := & -*_{\Fh} \partial_{\Fh} *_{\Fh} \; .
\end{eqnarray*}
Then $(d^*_{\Fh} \alpha) \tei \ell = d^*_\ell (\alpha \tei \ell)$ for each leaf $\ell$, where $d^*_{\ell}$ is the formal adjoint of $d_{\ell}$, and similarly for the other operators (cf. \cite{Wells}, Ch. V, Prop. 2.3, Prop. 2.4 and a footnote on p. 191). 

We form the Laplacians
\begin{eqnarray*}
  \Delta_{\Fh} & := & d_{\Fh} d^*_{\Fh} + d^*_{\Fh} d_{\Fh} \; , \\
\square_{\Fh} & := & \partial_{\Fh} \partial^*_{\Fh} + \partial^*_{\Fh} \partial_{\Fh} \; , \\
\oBox_{\Fh} & := & \opartial_{\Fh} \opartial^*_{\Fh} + \opartial^*_{\Fh} \opartial_{\Fh} \; .
\end{eqnarray*}
Finally, we define the operators $L_{\Fh}$ and $L^*_{\Fh}$ on $\Ah^{\hullet}_{\Fh} (X)$ by
\begin{eqnarray*}
  L_{\Fh} (\alpha) & := &\alpha \wedge \omega_{\Fh} \; , \\
L^*_{\Fh} & := & *^{-1}_{\Fh} L_{\Fh} *_{\Fh} \; .
\end{eqnarray*}
Observe that $L_{\Fh}$ and $L^*_{\Fh}$ are pointwise operators and that $L^*_{\Fh}$ is the pointwise adjoint of $L_{\Fh}$.

A form $\alpha \in \Ah^{\hullet}_{\Fh} (X)$ is called \textit{harmonic} if $\Delta_{\Fh} (\alpha) = 0$ and \textit{primitive} if $L^*_{\Fh} (\alpha) = 0$.

\begin{defn}
  A foliation $\Fh$ is called a \textit{K\"ahler foliation} if it is endowed with a complex structure $J$ and an Hermitian metric $h = S - 2i \omega_{\Fh}$ on $T^c \Fh$ such that
\[
d_{\Fh} \omega_{\Fh} = 0 \; .
\]
The leaves of a K\"ahler foliation are K\"ahler manifolds, and it is clear that the K\"ahler identities (\cite{Weil}, Ch. II, Sections 5 and 6, or \cite{Wells}, Ch. V, p. 191--195) hold in the foliated context. In particular,
\[
\Delta_{\Fh} = 2 \square_{\Fh} = 2 \oBox_{\Fh}
\]
is a real operator of bidegree $(0,0)$ which commutes with $L_{\Fh}$ and $L^*_{\Fh}$. Denote by $\Hh^r_{\Fh} (X)$ [resp. $\Hh^{p,q}_{\Fh} (X)$] the harmonic forms in $\Ah^r_{\Fh} (X)$ [resp. $\Ah^{p,q}_{\Fh} (X)$] and by $\Hh^r_{\Fh} (X)_0$ [resp. $\Hh^{p,q}_{\Fh} (X)_0$] the primitive elements in $\Hh^r_{\Fh} (X)$ [resp. $\Hh^{p,q}_{\Fh} (X)$]. Then as topological vector spaces:
\[
\Hh^{p,q}_{\Fh} (X) \cong \Hh^{q,p}_{\Fh} (X) \cong \Hh^{g-q, g-p}_{\Fh} (X) \cong \Hh^{g-p, g-q}_{\Fh} (X) \; .
\]
\end{defn}

From the Lefschetz decomposition theorem (\cite{Weil}, Th. 3 on p. 26 or \cite{Wells}, Th. 3.12 on p. 182) we obtain:
\begin{punkt}
  \label{t21} \rm
For a K\"ahler foliation $\Fh$ of dimension $2g$, we have
\[
\Hh^r_{\Fh} (X) = \bigoplus_{s \ge \max (r-g, 0)} L^s_{\Fh} (\Hh^{r-2s}_{\Fh} (X)_0)
\]
and, more precisely,
\[
\Hh^{p,q}_{\Fh} (X) = \bigoplus_{s \ge \max (p+q -g , 0)} L^s_{\Fh} (\Hh^{p-s,q-s}_{\Fh} (X)_0) \; .
\]
\end{punkt}

Since $L^{g-r}_{\Fh}$ is injective on $\Ah^r_{\Fh} (X)$ for $r < g$ (\cite{Wells}, Th. 3.12.c) on p. 182), \ref{t21} yields the following version of the hard Lefschetz theorem:

\begin{punkt}
  \label{t22} \rm
Let $\Fh$ be a K\"ahler foliation on $X$ of dimension $2g$. Then
\[
L^{g-r}_{\Fh} : \Hh^r_{\Fh} (X) \longrightarrow \Hh^{2g-r}_{\Fh} (X)
\]
is an isomorphism for $r < g$.
\end{punkt}

%%\newpage
%\newpage
%\input{sec3}
\section{K\"ahler foliations on closed manifolds}
\label{sec:3}
\begin{punkt}
  \label{t31} \rm
From now on, we consider foliations on closed manifolds. To fix some notations, we consider for the moment an arbitrary oriented and transversally oriented foliation $\Fh$ on the closed Riemannian manifold $X$. We then have star operators $*_{\Fh}$ and $*_{\perp}$ on the subbundles $\Lambda T^* \Fh$ and $\Lambda (T\Fh)^{\perp *}$ of $\Lambda T^*X$. On the manifold $X$ we take the orientation for which the volume form is given by
\[
\vol = *_{\perp} (1) \wedge *_{\Fh} (1) \; .
\]
We get a scalar product on $\Ah^r_{\Fh} (X) = \Gamma (X , \Lambda^r T^* \Fh \otimes \C)$ by
\[
(\alpha , \beta) := \int_X \langle \alpha , \beta \rangle \vol = \int_X \alpha \wedge * \obeta \; .
\]
Here $*$ is the star operator of $X$ which is equal to $*_{\perp} \otimes *_{\Fh}$ up to a sign. This sign is computed in \cite{AK}, Lemma 3.2, and is such that
\[
(\alpha , \beta) = \int_X *_{\perp} (1) \wedge \alpha \wedge *_{\Fh} (\obeta) \; .
\]
\end{punkt}

\begin{punkt}
  \label{t32} \rm
Now we specialize the situation of \ref{t31} by taking $\Fh$ to be a K\"ahler foliation on the closed orientable manifold $\Fh$. Let us retain the notations of section \ref{sec:2}. We can extend the Riemannian metric $S$ on $T\Fh$ to a Riemannian metric on $TX$ and choose an orientation for $T\Fh^{\perp}$. The last equation in \ref{t31} shows that
\[
\Ah^r_{\Fh} (X) = \bigoplus_{p+q=r} \Ah^{p,q}_{\Fh} (X)
\]
is an orthogonal decomposition.
\end{punkt}

\begin{punkt}
  \label{t33} \rm 
As in the non-foliated context, we consider a variant $Q$ of the scalar product $( \, , \,)$ on harmonic forms, as follows:\\
Given $\xi , \eta \in \Hh^r_{\Fh} (X)$, write $\xi = \sum_s L^s_{\Fh} \xi_s$ and $\eta = \sum_s L^s_{\Fh} \eta_s$ with primitive forms $\xi_s , \eta_s$ according to \ref{t21} and define
\[
Q (\xi , \eta) := \sum_s (-1)^{s + r (r+1) / 2} \int_X *_{\perp} (1) \wedge L^{g-r+2s}_{\Fh} (\xi_s \wedge \eta_s) \; .
\]
Obviously, we have
\[
Q (\xi , \eta) = (-1)^r Q (\eta , \xi) \; .
\]
Denote by $J$ the automorphism of $\Ah^{\hullet}_{\Fh} (X)$ which is multiplication by $i^{p-q}$ on $\Ah^{p,q}_{\Fh} (X)$. The proof of \cite{Wells}, Ch. V, Th. 6.1 shows that there are positive constants $c_s$ such that
\[
Q (\xi , J \oeta) = \sum_s c_s (L^s \xi_s , L^s \eta_s) \; .
\]
In particular, we have $Q (\xi , J\overline{\xi}) > 0$ for $\xi \neq 0$.

Moreover, $Q (J \xi , J \eta) = Q (\xi , \eta)$. 
\end{punkt}

\begin{punkt}
  \label{t34} \rm
As before, let $\Fh$ be a K\"ahler foliation of complex dimension $g$ on the closed manifold $X$; we assume now that $g$ is even. The bilinear form
\[
I (\alpha , \beta) = \int_X \alpha \wedge \beta
\]
on $\Hh^g_{\Fh} (X)$ is symmetric and real. Assume that $\dim \Hh^g_{\Fh} (X) < \infty$. We denote the signature of $I$ by sign $\Fh$ and call it the \textit{signature} of $\Fh$. By \ref{t22}, $\Hh^{\mathrm{even}}_{\Fh}$ is then finite dimensional, and the argument on p. 76--77 of \cite{Weil} shows that we have the following version of the Hodge index theorem:
\[
\sign (\Fh) = \sum_{p \equiv q (2)} (-1)^p \dim \Hh^{p,q}_{\Fh} (X) \; .
\]
\end{punkt}
%%\newpage
%\input{sec4}
\section{K\"ahler--Riemann foliations}
\label{sec:4}

So far, cohomology has not entered our discussion. To have a K\"ahler theory on the cohomology of a foliation, we have to restrict our foliations further: Let $(X, \Fh)$ be any foliation. Recall that a Riemannian metric on $X$ is called \textit{bundle-like} for $\Fh$, if the normal plane field of $\Fh$ is totally geodesic. Equivalently, the metric is bundle-like if for each pair of normal vector fields $Y$ and $Z$ which are infinitesimal automorphisms of $\Fh$ and for each tangential vector field $W$ we have $W \langle Y , Z \rangle = 0$. A foliation which admits a bundle-like metric is called \textit{Riemannian}.

Given an oriented and transversally oriented foliation $\Fh$ on the closed manifold $X$ and a bundle-like metric for $\Fh$, our leafwise adjoint $d^*_{\Fh}$ is the formal adjoint of the operator $d_{\Fh}$ on $X$ (\cite{AK}, Section 3), and a fundamental result of $\cite{AK}$ says that there is an orthogonal Hodge decomposition
\begin{equation}
  \label{eq:41}
  \Ah^{\hullet}_{\Fh} (X) = \Hh^{\hullet}_{\Fh} (X) \oplus \overline{\imm d_{\Fh}} \oplus \overline{\imm d^*_{\Fh}} \; .
\end{equation}
Here, the bar denotes closure in the Fr\'echet space $\Ah^{\hullet}_{\Fh} (X)$. In particular, the reduced leafwise complex cohomology
\[
\oH^{\hullet}_{\Fh} (X , \C) = \ker d_{\Fh} / \overline{\imm d_{\Fh}}
\]
is isomorphic to $\Hh^{\hullet}_{\Fh} (X)$. This Hodge decomposition does not hold in general for non-Riemannian foliations \cite{DS}.

\begin{remark}
  Denoting the dimension of $\Fh$ by $m$, we get a linear form $\nu$ on $\Gamma (X , \Lambda^m T^* \Fh)$ by
\[
\nu (\omega) := \int *_{\perp} (1) \wedge \omega = (\omega , *_{\Fh} (1)) \; .
\]
Since $d^*_{\Fh}$ is the formal adjoint of $d_{\Fh}$, we see that $\nu (d_{\Fh} \beta) = 0$; hence $\nu$ induces a (continuous) linear form on $H^m_{\Fh} (X)$, the leafwise real $m$-dimensional cohomology group. This functional in turn induces a linear form
\[
\nu : \oH^m_{\Fh} (X) \longrightarrow \R
\]
on the corresponding reduced group. We get a scalar product on every $\oH^n_{\Fh} (X)$ by setting:
\begin{eqnarray*}
  (a,b) & = & \nu (a \cup *_{\Fh} b) \\
& = & \int_X \langle \Hh (a) , \Hh (b) \rangle_{\Fh} \vol \; .
\end{eqnarray*}
Here $\Hh (a) , \Hh (b)$ are the harmonic representatives of $a$ and $b$ and $*_{\Fh} a := [*_{\Fh} \Hh (a)]$. In fact under the isomorphism
\[
\oH^n_{\Fh} (X) \cong \Hh^n_{\Fh} (X)
\]
this scalar product becomes the restriction of the scalar product on forms from \ref{t31}.
\end{remark}

The bilinear form
\[
(a,b) \longmapsto \nu (a \cup b)
\]
on $\oH^{\hullet}_{\Fh} (X)$ is obviously non-degenerate in the sense that for every $a \in \oH^{\hullet}_{\Fh} (X)$ with $a \neq 0$, there exists $b$ with $\nu (a \cup b) \neq 0$. If $\Fh$ has a dense leaf, we have $\oH^m_{\Fh} (X) \cong \oH^0_{\Fh} (X) \cong \R$, and there is the following version of Poincar\'e duality:

\begin{punkt}
  \label{t42} \rm
Let $\Fh$ be an orientable and transversely orientable Riemannian foliation of dimension $m$ on the closed manifold $X$. Assume that $\Fh$ has a dense leaf. For any isomorphism $\int : \oH^m_{\Fh} (X) \to \R$, the bilinear form $(a,b) \mapsto \int a \cup b$ is non-degenerate.
\end{punkt}

\begin{defn}
  A K\"ahler foliation $(X , \Fh , h = S - 2 i \omega_{\Fh})$ which is also a Riemannian foliation is called a \textit{K\"ahler--Riemann foliation}. A Riemannian metric on $X$ is \textit{K\"ahler-bundle-like} if it is bundle-like and if its restriction to $T \Fh$ agrees with $S$.
\end{defn}

A K\"ahler--Riemann foliation always admits a K\"ahler-bundle-like metric. Hence, given a K\"ahler--Riemann foliation $\Fh$ on the closed orientable manifold $X$, we obtain a decomposition
\begin{equation}
  \label{eq:43}
  \oH^r_{\Fh} (X , \C) = \bigoplus_{p+q = r} H^{p,q}
\end{equation}
where $H^{p,q}$ is the image of $\Hh^{p,q}_{\Fh} (X)$ under the natural isomorphism. The subspace $H^{p,q}$ is closed in $\oH^r_{\Fh} (X,\C)$, and we have $\overline{H^{p,q}} = H^{q,p}$; in this equation, the bar denotes complex conjugation.

In particular, a necessary condition for the existence of a K\"ahler structure on a Riemannian foliation is this: The reduced leafwise cohomology groups in odd degree must be of even dimension if they are finite dimensional.

In particular for Riemannian foliations by (real) surfaces $\dim \oH^1_{\Fh} (X)$ is even if finite since such foliations always carry a K\"ahler structure: A metric on $T\Fh$ determines an almost complex structure on $T\Fh$ which is integrable since the leaves are $2$-dimensional. The associated K\"ahler form is closed for degree reasons.

\begin{prop}
  \label{t44}
$H^{p,q}$ consists of those reduced cohomology classes which can be represented by a cocycle contained in $\Ah^{p,q}_{\Fh} (X)$. In particular, the decomposition (2) only depends on the complex structure $J$ and not on the K\"ahler form.
\end{prop}

\begin{proof}
  Let $\pi : \Ah^{\hullet}_{\Fh} (X) \to \Hh^{\hullet}_{\Fh} (X)$ be the projection according to (1). Since the decompositions (1) and $\Ah^{\hullet} = \bigoplus_{p,q} \Ah^{p,q}$ are orthogonal, it is clear that $\pi (\Ah^{pq}_{\Fh} (X)) = \Hh^{p,q}_{\Fh} (X)$.

Therefore, given an element $\beta \in \Ah^{p,q}_{\Fh} (X)$ with $d_{\Fh} \beta = 0$, the cohomology class of $\beta$ is also represented by $\pi (\beta)$, hence is contained in $H^{p,q}$. The rest of \ref{t44} is obvious.
\end{proof}

For a K\"ahler--Riemann foliation of dimension $2g$ on a closed orientable manifold $X$, the results of section \ref{sec:2} and \ref{sec:3} can be translated from forms to cohomology:

Forming the cup product with the reduced K\"ahler class $[\omega_{\Fh}] \in \oH^2_{\Fh} (X)$ defines a real operator $L : \oH^{\hullet}_{\Fh} (X , \C) \to \oH^{\hullet}_{\Fh} (X , \C)$ of bidegree $(1,1)$, and \ref{t22} yields a hard Lefschetz theorem:

\begin{prop}
  \label{t45} \rm
For $r < g$, we have an isomorphism
\[
L^{g-r} : \oH^r_{\Fh} (X , \C) \longrightarrow \oH^{2g-r}_{\Fh} (X , \C) \; .
\]
In particular, $[\omega_{\Fh}]^g \neq 0$.
\end{prop}

Also, \ref{t34} gives a Hodge index theorem for reduced leafwise cohomology.

We call a class $a \in \oH^r_{\Fh} ( X , \C)$ \textit{primitive} if it corresponds to a primitive harmonic form $\alpha \in \Hh^r_{\Fh} (X)_0$ under the natural isomorphism $\Hh^r_{\Fh} (X) \cong \oH^r_{\Fh} (X , \C)$. For $r \le g$, the class $a$ is primitive iff $L^{g-r+1} a = 0$; for $r > g$ there are no non-zero primitive classes. 

Denoting the set of primitive $r$-classes by $P^r_{\Fh} (X)$, we have by \ref{t21}:
\begin{equation}
  \label{eq:46}
  \oH^r_{\Fh} (X , \C) = \bigoplus_s L^s (P^{r-2s}_{\Fh} (X)) \; .
\end{equation}
The bilinear form $Q$ introduced in \ref{t33} yields a bilinear form $Q$ on $\oH^r_{\Fh} (X , \C)$ which can be described as follows: Given $a,b \in \oH^r_{\Fh} (X , \C)$, write $a = \sum L^s a_s , b = \sum L^s b_s$ according to (3). Then we have
\[
Q (a,b) = \sum_s (-1)^{s + r (r+1)/2} \nu (L^{g-r+2s} (a_s \cup b_s))
\]
where $\nu$ is as in the discussion preceding \ref{t42}.

If $\Fh$ has a dense leaf, the non-zero functional $\nu$ on $\oH^{2g}_{\Fh} (X , \C)$ is determined up to a constant, and we see that in this case $Q$ is essentially determined by the K\"ahler class $[\omega_{\Fh}]$. 

The preceeding discussion also holds for the reduced leafwise cohomology with real coefficients.

Let $J$ be the automorphism of $\oH^n_{\Fh} (X , \C)$ which is multiplication by $i^{p-q}$ on $H^{pq}$. It follows from \ref{t33} that 
\begin{equation}
  \label{eq:47}
  T (a,b) = Q (a, J\bar{b})
\end{equation}
defines a scalar product on $\oH^n_{\Fh} (X , \C)$.

From \ref{t33} we also conclude:

\begin{punkt}
  \label{t48} \rm
$Q$ defines a polarization of weight $r$ of the real Hodge structure $P^r_{\Fh} (X)$.

Here, we define real Hodge structures and their polarizations as Deligne (\cite{De}, (2.1.5), (2.1.16)) except that we do not require that our vector spaces are finite dimensional.
\end{punkt}

\begin{punkt}
  \label{t411} \rm
In this subsection we generalize a theorem of Serre \cite{S} on a K\"ahler analogue of the Weil conjectures to our context. So assume that $\Fh$ is a K\"ahler--Riemann foliation of complex dimension $g$ with a dense leaf on the oriented manifold $X$.
\end{punkt}

\begin{prop}
  \label{t412}
Let $f : X\to X$ be a smooth map which maps leaves of $\Fh$ holomorphically into leaves. Assume that the induced map $f^*$ on $\oH^{\hullet}_{\Fh} (X , \C)$ satisfies the relation:
\[
f^* [\omega_{\Fh}] = q [\omega_{\Fh}]
\]
for some positive real number $q > 0$. Then on $\oH^n_{\Fh} (X, \C)$ we have:
\[
f^* = q^{n/2} U_n
\]
where $U_n$ is a unitary endomorphism of $\oH^n_{\Fh} (X , \C)$ with respect to the scalar product $T$ on cohomology defined by (4). In particular the spectral values of $f^*$ have absolute value $q^{n/2}$.
\end{prop}

\begin{proof}
  Define an endomorphism $U$ of the vector space $\oH^{\hullet}_{\Fh} (X , \C)$ by setting $U_n = q^{-n/2} f^*$ in degree $n$. Then $U$ is an algebra endomorphism of $\oH^{\hullet}_{\Fh} (X , \C)$. We know that $[\omega_{\Fh}]^g \neq 0$ by \ref{t45}. Moreover since $\Fh$ has a dense leaf, $\oH^{2g}_{\Fh} (X , \C)$ is one dimensional. Hence $U_{2g} = \id$. Since $U$ commutes with the Lefschetz operator $L$, we find
\[
Q (Ua , Ub) = Q (a,b) \; .
\]
Since $f^*$ respects the real structure $\oH^{\hullet}_{\Fh} (X)$ of $\oH^{\hullet}_{\Fh} (X , \C)$ and since $f^* H^{pq} \subset H^{pq}$ it follows that $U$ commutes with $J$ and with complex conjugation. Thus
\[
T (Ua , Ub) = T (a,b) \; .
\]
\end{proof}

\begin{punkt}
  \label{t49} \rm
We now connect the spaces $H^{pq}$ in $\oH^n_{\Fh} (X, \C)$ to the reduced leafwise Dolbeault cohomology groups and make a number of further comments. So let $\Fh$ be as before a K\"ahler--Riemann foliation on the oriented manifold $X$ and fix a K\"ahler bundle-like metric. Using the equation $\Delta_{\Fh} = 2\oBox_{\Fh}$ and passing to $(p,q)$-types in the Hodge decomposition of \cite{AK}:
\[
\Ah^{\hullet}_{\Fh} (X) = \Ker \Delta^{\hullet}_{\Fh} \oplus \overline{\imm \Delta^{\hullet}_{\Fh}} \; ,
\]
we obtain the orthogonal decomposition:
\[
\Ah^{pq}_{\Fh} (X) = \Ker \oBox^{pq}_{\Fh} \oplus \overline{\imm \oBox^{pq}_{\Fh}} \; .
\]
Since $\imm \opartial_{\Fh}$ and $\imm \opartial^*_{\Fh}$ are orthogonal, we get the following Hodge decomposition:
\[
\Ah^{pq}_{\Fh} (X) = \Hh^{pq}_{\Fh} (X) \oplus \overline{\imm \opartial^{p,q-1}_{\Fh}} \oplus \overline{\imm (\opartial^{p,q+1}_{\Fh})^*} \; .
\]
In particular the reduced leafwise Dolbeault cohomology $\oH^q (X , \Omega^p_{\Fh})$ is isomorphic to $\Hh^{pq}_{\Fh} (X)$ and hence to the subspace $H^{pq}$ of $\oH^{p+q}_{\Fh} (X, \C)$.

The Hodge $\overline{*}_{\Fh}$-operator induces an anti-isomorphism:
\[
\oH^q (X , \Omega^p_{\Fh}) \silo \oH^{g-q} (X , \Omega^{g-p}_{\Fh}) \; .
\]
Moreover, it follows in the case of a dense leaf that the $\cup$-product pairing:
\[
\oH^q (X , \Omega^p_{\Fh}) \times \oH^{g-q} (X, \Omega^{g-p}_{\Fh}) \overset{\cup}{\longrightarrow} \oH^g (X , \Omega^g_{\Fh}) \cong H^{gg} = \oH^{2g}_{\Fh} (X) \overset{\nu}{\silo} \R
\]
is non-degenerate, a weak version of Serre duality.

Define the Hodge filtration $F^i$ on $\oH^n_{\Fh} (X , \C)$ by:
\[
F^i = \imm (H^n (X , \Omega^{\ge i}_{\Fh}) \longrightarrow \oH^n (X , \Omega^{\hullet}_{\Fh}) = \oH^n_{\Fh} (X , \C)) \; .
\]
Note here that the complex $\Omega^{\hullet}_{\Fh}$ is quasi-isomorphic to $\Ch_{\Fh}$. 
\end{punkt}

\begin{prop}
  \label{t410}
The subspace $F^i \subset \oH^n_{\Fh} (X , \C)$ is closed and we have
\[
F^i = \bigoplus_{p \ge i} H^{pq} \; .
\]
\end{prop}

\begin{proof}
  The diagram
\[
\begin{CD}
  \Omega^{\ge i}_{\Fh} @>>> \Omega^{\hullet}_{\Fh} \\
@VVV @VVV \\
s \Ah^{\ge i, \hullet}_{\Fh} @>>> s \Ah^{\hullet\hullet}_{\Fh}
\end{CD}
\]
is commutative and the vertical arrows are quasi-isomorphisms as follows from \ref{t14}. Here $s$ denotes the associated simple complex of a double complex. Thus $F^i$ can be identified with the image of the natural map
\[
H^n \Big( \bigoplus_{p + q = \hullet \atop p \ge i} \Ah^{pq}_{\Fh} (X) , d_{\Fh} \Big) \longrightarrow \oH^n \Big( \bigoplus_{p+q = \hullet} \Ah^{pq}_{\Fh} (X) , d_{\Fh} \Big) \; .
\]
Its image consists of those cohomology classes that can be represented by closed forms of the type:
\[
\sum_{p+q = n \atop p \ge i} \alpha_{pq} \quad \mbox{with} \; \alpha_{pq} \in \Ah^{pq}_{\Fh} (X) \; .
\]
Passing to the harmonic projections of the $\alpha_{pq}$, we see that every form in the image can be represented by a form
\[
\sum_{p+q = n \atop p \ge i} h_{pq} \quad \mbox{with} \; h_{pq} \in\Ah^{pq}_{\Fh} (X) \; \mbox{and} \; dh_{pq} = 0 \; .
\]
Hence 
\[
F^i = \bigoplus_{p \ge i} H^{pq}
\]
and in particular $F^i$ is a closed subspace, since the $H^{pq}$'s are. 
\end{proof}

%\input{appendix}
%%\newpage
%\input{sec5}
\section{The Lie algebra attached to a K\"ahler--\\
Riemann foliation}
\label{sec:5}

There is a construction due to Looijenga and Lunts \cite{LL} attaching a Lie algebra to each compact K\"ahler manifold. We extend this construction to K\"ahler--Riemann foliations; new phenomena can occur because the cohomology may be of infinite dimension.

In this section we consider only complex Lie algebras.

Let $\Fh$ be a $2g$-dimensional K\"ahler--Riemann foliation of the closed orientable manifold $X$. Let $B$ be the endomorphism of the graded vector space $\oH^{\hullet} := \oH^{\hullet}_{\Fh} (X , \C)$ which is multiplication by $g-r$ on $\oH^r$. For $\omega \in \oH^2$, let $L_{\omega}$ be the endomorphism of $\oH^{\hullet}$ of degree $2$ given by $L_{\omega} (x) = x \cup \omega$. By $\Kh_{\Fh} (X)$ we denote the non-empty set of all those $\omega \in \oH^2$ for which $L^k_{\omega}$ is an isomorphism of $\oH^{g-k}$ onto $\oH^{g+k}$ for all $k \ge 0$.

\begin{punkt}
  \label{t51} \rm
For $\omega \in \Kh_{\Fh} (X)$, there exists a uniquely determined endomorphism $\Lambda_{\omega}$ of $\oH^{\hullet}$ of degree $-2$ such that
\[
[\Lambda_{\omega} , L_{\omega}] = B \; .
\]
\end{punkt}

\begin{proof}
  As observed in \cite{LL}, the assertion follows from the Jacobson--Morozov lemma in the finite dimensional case. The general case can be easily deduced.
\end{proof}

The two elements $L_{\omega}$ and $\Lambda_{\omega}$ generate a Lie subalgebra of $\End (\oH^{\hullet})$ isomorphic to $\esl_2$. 

Let $\eg = \eg (\Fh)$ be the Lie algebra generated by all the $L_{\omega}$ and $\Lambda_{\omega}$ with $\omega \in \Kh_{\Fh} (X)$. If $\dim \oH^{\hullet} < \infty$ then $\oH^{\hullet}$ is an $\oH^2$-Lefschetz-module in the sense of \cite{LL}, that is, $\eg$ is semisimple.

Of course, one can consider the even and the odd parts of $\oH^{\hullet}$ separately: Let $\Kh_{\mathrm{even}}$ be the set of all $\omega \in \oH^2$ with the following property: For all $k$ with $k \equiv g (\mod 2) , L^k_{\omega}$ is an isomorphism of $\oH^{g-k}$ onto $\oH^{g+k}$. For $\omega \in \Kh_{\mathrm{even}}$, there is an associated endomorphism $\Lambda_{\omega}$ of $\oH^{\mathrm{even}}$. Let $\eg_{\mathrm{even}}$ be the Lie subalgebra of $\End (\oH^{\mathrm{even}})$ generated by the $L_{\omega}$ and $\Lambda_{\omega}$ with $\omega \in \Kh_{\mathrm{even}}$. Similarly, one can define $\eg_{\mathrm{odd}}$.

As an example, let us determine the Lie algebra $\eg_{\even}$ in the case that $g = 2$ and that $\Fh$ has a dense leaf. In this situation, $\oH^0 = \C$ and we can choose an isomorphism $u \mapsto \int u$ of $\oH^4$ onto $\C$. We obtain a symmetric bilinear form $\varphi$ on $\oH^{\even}$ by
\begin{eqnarray*}
  \varphi \tei \oH^k \times \oH^h = 0 & \mbox{if} & k + h \neq 4 \; , \\
\varphi (u,v) = \int uv & \mbox{if} & u \in \oH^4 , v \in \oH^0 \; , \\
\varphi (u,v) =- \int uv & \mbox{if} & u, v \in \oH^2 \; .
\end{eqnarray*}
Let $\so (\oH^{\even} , \varphi)$ be the Lie algebra of all endomorphisms $A$ of $\oH^{\even}$ with
\[
\varphi (Au , v) + \varphi (u, Av) = 0
\]
for all $u,v \in \oH^{\even}$, and let $\so_{\fin} (\oH^{\even} , \varphi)$ be the Lie subalgebra of $\so (\oH^{\even} , \varphi)$ consisting of all endomorphisms of finite rank.

\begin{theorem}
  \label{t52}
If $g = 2$ and if $\Fh$ has a dense leaf, we have
\[
\eg_{\even} = \so_{\fin} (\oH^{\even} , \varphi) \; .
\]
If $\dim \oH^{\even} = \infty$, this is a simple Lie algebra.
\end{theorem}

\begin{proof}
  Let us abbreviate $\so_{\fin} (\oH^{\even} , \varphi)$ by $\eg$. In our situation, 
\[
\Kh_{\even} = \{ \omega \in \oH^2 \tei \omega^2 \neq 0 \} \; .
\]
For $\omega \in \Kh_{\even}$, the homomorphism $\Lambda_{\omega} : \oH^4 \to \oH^2$ is given by $\Lambda_{\omega} (\omega^2) = 2 \omega$, and
\[
\Lambda_{\omega} : \oH^2 = \C \cdot \omega \oplus \ker L_{\omega} \longrightarrow \oH^0
\]
is given by
\[
\begin{array}{l}
\Lambda_{\omega} (\omega) = 2 \; , \\
\ker (\Lambda_{\omega} \tei \oH^2) = \ker (L_{\omega} \tei \oH^2) \; .
\end{array}
\]
Since $L_{\omega}$ and $\Lambda_{\omega}$ have rank $2$ on $\oH^{\even}$ and since obviously 
\[
L_{\omega} , \Lambda_{\omega} \in \so (\oH^{\even} , \varphi)\; ,
\]
we see that $\eg_{\even} \subseteq \eg$. Now observe the following three points, \\
using \ref{t42} for (a):\\
(a) Given a finite dimensional subspace $U$ of $\oH^{\hullet}$, there exists a finite dimensional subspace $W$ of $\oH^{\hullet}$ containing $U$ such that the restriction $\varphi_W$ of $\varphi$ to $W \times W$ is non-degenerate. \\
(b) If $W$ is a finite dimensional subspace of $\oH^{\hullet}$ such that $\varphi_W$ is non-degenerate, there is a decomposition $\oH^{\hullet} = W \oplus W^{\perp}$ where $W^{\perp}$ is formed with respect to $\varphi$. Let
\[
\eg_W := \{ A \in \eg \tei A (W) \subseteq W \; \mbox{and} \; A \tei W^{\perp} = 0 \} \; .
\]
Then $\eg_W$ is a Lie subalgebra of $\eg$ which is isomorphic to $\so (W , \varphi_W)$. \\
(c) For $A \in \eg$ and a finite dimensional subspace $W$ of $\oH^{\hullet}$ containing $\imm A$ and such that $\varphi_W$ is non-degenerate, we have $A \in \eg_W$.

To see that $\eg \subseteq \eg_{\even}$, consider an element $A \in \eg$ and choose a finite dimensional subspace $W$ such that $\varphi_W$ is non-degenerate and $A \in \eg_W$. We may assume that $W$ is of the form $W = \oH^0 \oplus W^2 \oplus \oH^4$ with $W^2 \subseteq \oH^2$. Then $W$ is a graded algebra, and it is clear (cf. \cite{LL}, the beginning of section 4) that the Lie algebra attached to it is $\so (W , \varphi_W)$. Hence $A \in \eg_{\even}$.

To see that $\eg$ is simple if $\dim \oH^{\hullet} = \infty$, assume that there is a proper ideal $\ea$ of $\eg$. Take elements $A \in \ea , A \neq 0$, and $C \in \eg , C \notin \ea$. Let $W$ be a finite dimensional subspace of $\oH^{\hullet}$ containing $\imm A$ and $\imm C$ such that $\varphi_W$ is non-degenerate. Then $A, C \in \eg_W$. Hence $\ea \cap \eg_W$ is a proper ideal in $\eg_W$. Assuming that $\dim W > 4$, the Lie algebra $\eg_W$ is simple. This is a contradiction.
\end{proof}

\begin{remarks}
  \label{t53}
\rm If there is no dense leaf, the situation looks quite different: For example, let $X = S^1 \times S^2$ with leaves $\{ a \} \times S^2$. Then 
\[
\eg (\Fh) = \esl_2 (C^{\infty} (S^1)) \; ,
\]
considered as a complex Lie algebra. It has very many ideals.

Even if there is a dense leaf and if $g \ge 3$, it is not likely that $\eg (\Fh)$ is always a finite product of simple Lie algebras.
\end{remarks}

\section{Examples of K\"ahler--Riemann foliations}
\label{sec:6}

In this section we consider a class of examples that were suggested to us by E. Ghys. We begin with the following general setup which will be specialized to our needs later.

Let $G$ be a real Lie group with a discrete subgroup $\Gamma$ and a compact connected subgroup $K$ such that $\Gamma$ operates freely on the manifold $S = G / K$. Then $X = \Gamma \ohne S$ is a manifold as well. Consider a connected normal sub Lie group $G_1$ of $G$. Then $K_1 = G_1 \cap K$ is a compact normal subgroup of $K$.

The orbits of the left $G_1$-operation on $S$ define a foliation $\Fh_S$ on $S$ of dimension $\dim G_1 / K_1$. This follows because the isotropy groups 
\[
(G_1)_{gK} = gK g^{-1} \cap G_1 = g K_1 g^{-1}
\]
all have the same dimension, $G_1$ being normal in $G$, \cite{Go} 1.12 iii). The $\Gamma$-operation on $S$ carries leaves of $\Fh_S$ into leaves:
\[
\gamma (G_1 g K) = G_1 \gamma g K \; ,
\]
again since $G_1$ is normal. Thus $\Fh_S$ descends to a foliation $\Fh$ on $X$ with leaves
\[
\Gamma G_1 gK = \Gamma g G_1 K \; .
\]
The space of leaves on $X$ is therefore:
\begin{equation}
  \label{eq:5}
  X / \Fh = \Gamma \ohne G / G_1 K \; .
\end{equation}
There are natural isomorphisms of vector bundles on $S$:
\begin{equation}
  \label{eq:6}
  G \times_K \eg / \ek \silo TS \quad \mbox{and} \quad G \times_K \eg_1 / \ek_1 \silo T\Fh_S \; ,
\end{equation}
where $K$ acts via $\Ad$ on $\eg / \ek$ and $\eg_1 / \ek_1$. They are given by mapping $[g , \ou]$ to $(T_{\ove} L_g) (\ou)$ where $L_g$ denotes left multiplication with $g$ on $S$ and $\ou \in \eg / \ek = T_{\ove} (G / K)$ resp.
\[
\ou \in \eg_1 / \ek_1 = T_{\ove} (G_1 / K_1) = T_{\ove} (G_1 K / K) = T_{\ove} \Fh_S \; .
\]
The isomorphisms (6) are left $G$-equivariant and hence induce isomorphisms
\begin{equation}
  \label{eq:7}
  (\Gamma \ohne G) \times_K \eg / \ek \silo TX \quad \mbox{and} \quad (\Gamma \ohne G) \times_K \eg_1 / \ek_1 \silo T\Fh \; .
\end{equation}
From (7) we get isomorphisms
\[
\Gamma (X , TX) \silo C^{\infty} (\Gamma \ohne G , \eg / \ek)^K \quad \mbox{and} \quad \Gamma (X , T \Fh) \silo C^{\infty} (\Gamma  \ohne G , \eg_1 / \ek_1)^K \; .
\]
If $f_A$ is the $K$-invariant function corresponding to the vector field $A$ then we have explicitely:
\[
f_A (\Gamma g k^{-1}) = \Ad (k) f_A (\Gamma g) \quad \mbox{for all} \; k \in K , g \in G \; .
\]
The map
\[
C^{\infty} (\Gamma \ohne G , \eg)^K \longrightarrow C^{\infty} (\Gamma \ohne G , \eg / \ek)^K
\]
has a section since $K$ is compact. Hence $f_A$ can be lifted to a $K$-invariant function $\tilde{f}_A : \Gamma \ohne G \to \eg$. Under the isomorphism
\[
(\Gamma \ohne G) \times \eg \silo T (\Gamma \ohne G)
\]
the function $\tilde{f}_A$ corresponds to a vector field $\tilde{A}$ on $\Gamma \ohne G$. In these terms we have the formula
\begin{equation}
  \label{eq:8}
  f_{[A,B]} (\Gamma g) = \tilde{A} (\tilde{f}_B) (\Gamma g) - \tilde{B} (\tilde{f}_A) (\Gamma g) + [\tilde{f}_A (\Gamma g) , \tilde{f}_B (\Gamma g)] \mod \ek \; .
\end{equation}
Choose any $K$-invariant scalar product $(\, , \,)$ on $\eg / \ek$. Via (7) it induces a Riemannian metric $g$ on $TX$ such that:
\[
g (A, B)_{\Gamma gK} = (f_A (\Gamma g) , f_B (\Gamma g)) \; .
\]
The following assertion must be well known:

\begin{prop}
  \label{t61}
The Riemannian metric $g$ is bundle-like for $\Fh$. In particular $\Fh$ is Riemannian.
\end{prop}

\begin{proof}
  We have the $K$-invariant decomposition:
\[
\eg / \ek = (\eg_1 + \ek) / \ek \oplus V
\]
where $V$ is the orthogonal complement of $(\eg_1 + \ek) / \ek$. 

It induces a $K$-invariant decomposition:
\begin{equation}
  \label{eq:9}
  \eg = (\eg_1 + \ek) \oplus V 
\end{equation}
where we have identified $V$ with its image in $\eg / \ek$. 

A vector field $A$ on $X$ is orthogonal to $T \Fh$ if and only if $f_A$ takes values in $V$. In this case $\tilde{f}_A$ takes values in $V$ as well.

Now assume in addition that $A$ normalizes $T\Fh$. Then for every $B \in \Gamma (X , T\Fh)$ the function $f_{[A,B]}$ takes values in $\eg_1 + \ek / \ek$. Since $\tilde{f}_B$ takes values in $\eg_1$, by formula (8) this means that $\tilde{B} (\tilde{f}_A)$ takes values in $\eg_1 + \ek$. On the other hand, since $\tilde{f}_A$ is $V$-valued, $\tilde{B} (\tilde{f}_A)$ is $V$-valued as well, so that $\tilde{B} (\tilde{f}_A) = 0$ for all $\tilde{B}$. This means that $\tilde{f}_A$ is constant on the foliation of $\Gamma \ohne G$ by the cosets of $G_1 K$. Thus $\tilde{f}_A$ and hence $f_A$ is a $V$-valued function on $\Gamma \ohne G / G_1 K = X / \Fh$. 

Therefore if $A_1 , A_2$ are two vector fields orthogonal to $T\Fh$ and normalizing $T\Fh$ it follows that
\[
g (A_1 , A_2) = (f_A , f_B)
\]
is a basic function. Thus $g$ is bundle like.
\end{proof}

Let $C^{\hullet} (\eg , K , V)$ be the complex calculating $(\eg , K)$-cohomology of $V$ c.f. \cite{BW} I \S\,5.

Using (7) we have natural topological isomorphisms:
\begin{eqnarray*}
  \Ah^p_{\Fh} (X) & = & \Gamma (X , \Lambda^p T^* \Fh \otimes \C) \\
& = & C^{\infty} (\Gamma \ohne G , \Lambda^p (\eg_1 / \ek_1)^* \otimes \C) \\
& = & \Hom_K (\Lambda^p (\eg_1 / \ek_1) , C^{\infty} (\Gamma \ohne G)) \\
& = & C^p (\eg_1 , K_1 , C^{\infty} (\Gamma \ohne G))^{K / K_1} \; .
\end{eqnarray*}
One checks that we get a topological isomorphism of complexes
\[
\Ah^{\hullet}_{\Fh} (X) = C^{\hullet} (\eg_1 , K_1 , C^{\infty} (\Gamma \ohne G))^{K / K_1} \; .
\]

Since $K / K_1$ is a compact group, it follows from this that
\begin{equation}
  \label{eq:10}
  \oH^{\hullet}_{\Fh} (X , \C) = \oH^{\hullet} (\eg_1 , K_1 , C^{\infty} (\Gamma \ohne G))^{K / K _1} \; .
\end{equation}
If $\Gamma \ohne G$ is compact it carries a unique invariant probability measure. We have the decomposition
\[
L^2 (\Gamma \ohne G) = \hat{\bigoplus_{\pi \in {\hat{G}}}} H (\pi)
\]
of the ``right regular'' representation of $G$ on the complex $L^2$-space of $\Gamma \ohne G$ into isotypic components $H (\pi)$. Here $\pi$ runs over the unitary dual $\hat{G}$ of $G$. 

To state the next result let us introduce the following notation. For a topological vector space $H$ and linear subspaces $H_i \subset H$ the equation
\[
H = \bar{\bigoplus_i} H_i
\]
by definition means the following:\\
i) The algebraic direct sum $\bigoplus\limits_i H_i$ is a dense subspace of $H$.\\
ii) Every element of $H$ can be written as an unconditionally convergent series $\sum_i h_i$ with $h_i \in H_i$.

Let $H (\pi)^{\infty}$ be the subspace of smooth vectors in $H (\pi)$.

\begin{prop}
  \label{t62} If $\Gamma \ohne G$ is compact we have
\[
\oH^{\hullet} (\eg_1 , K_1 , C^{\infty} (\Gamma \ohne G)) = \bar{\bigoplus_{\pi \in {\hat{G}}}} \oH^{\hullet} (\eg_1 , K_1 , H (\pi)^{\infty}) \; .
\]
\end{prop}

\begin{proof}
  A sequence $(E_i)_{i \ge 1}$ of non-trivial linear subspaces of a Fr\'echet space $E$ is called a Schauder basis of subspaces of $E$ if the following conditions hold:\\
a) Every element $v \in E$ can be written uniquely as a convergent series $v = \sum^{\infty}_{i=1} v_i$ with $v_i \in E_i$.\\
b) The projections $E \to E_i , v \mapsto v_i$ resulting from a) are continuous.\\
If the series in a) converge unconditionally the basis $(E_i)$ is called unconditional.

The proof of the following auxiliary result is straightforeward and will be omitted.

\begin{lemma}
  \label{t621}
Let $C^{\hullet}$ be a complex of Fr\'echet spaces with continuous differentials. Consider closed subcomplexes $C^{\hullet}_i \subset C^{\hullet}$ for $i \ge 1$ such that for every $p$ the sequence $(C^p_i)_{i \ge 1}$ is an unconditional Schauder basis of subspaces of $C^p$. Then we have:
\[
\oH^p (C^{\hullet}) = \bar{\bigoplus_{i \ge 1}} \oH^p (C^{\hullet}_i) \; .
\]
\end{lemma}

Using the Sobolev embedding lemma one checks that the subspaces $H (\pi)^{\infty}$ of $C^{\infty} (\Gamma \ohne G) = L^2 (\Gamma \ohne G)^{\infty}$ with $H (\pi)^{\infty} \neq 0$ form an unconditional Schauder basis of subspaces of $C^{\infty} (\Gamma \ohne G)$.

Consequently the subcomplexes
\[
C^{\hullet} (\eg_1 , K_1 , H (\pi)^{\infty})^{K / K_1} 
\]
of $C^{\hullet} (\eg_1 , K_1 , C^{\infty} (\Gamma \ohne G))^{K / K_1}$ satisfy the assumptions of the lemma and the assertion of \ref{t62} follows.
\end{proof}

In order to get K\"ahler--Riemann foliations in the above setting we assume that $G_1 / K_1$ is a bounded symmetric space. If $G_1$ is semisimple then the condition that $G_1$ be normal in $G$ implies that $\eg$ decomposes into a product $\eg = \eg_1 \times \eg_2$. Hence we consider from now on the following setup:\\
$G = G_1 \times G_2$ \quad  where $G_1$ is a connected reductive
Lie group in the sense of \cite{BW} 0.3.1 with compact center. $G_2$ is any Lie group. \\
$K = K_1 \times K_2$ \quad where $K_1$ is a maximal compact subgroup of $G_1$ and $K_2$ is a compact subgroup of $G_2$.\\
Let $\Theta$ be the Cartan involution associated to $K_1$. We have the $K_1$-invariant Cartan decomposition
\[
\eg_1 = \ek_1 \oplus \ep_1 \quad \mbox{where} \; \ep_1 = \{ x \in \eg \tei \Theta (x) = -x \} \; .
\]
We assume that $G_1 / K_1$ carries an invariant complex structure. There is then an element $z_0$ in the center of $\ek_1$ such that $J = \ad z_0 \, |_{\ep_1}$ defines a complex structure on $\ep$ which is invariant under $K_1$. Hence $J$ induces an almost complex structure of $G_1 / K_1$ which agrees with the given invariant complex structure. The complexification of $\ep_1$ decomposes into the $\pm i$-eigenspaces %\linebreak
of $J$:
\[
\ep_{1 \C} = \ep^+_1 \oplus \ep^-_1 \; .
\]

Let $\Gamma$ be a cocompact lattice in $G$ which acts without fixed points on $G / K$. Then 
\[
X = \Gamma \ohne G /K = \Gamma \ohne G_1 \times G_2 / K_1 \times K_2
\]
is a compact manifold. The Riemannian foliation $\Fh$ is given by the images of $G_1 / K_1 \times \{ g_2 K \}$ for $g_2 K \in G_2 / K_2$ under the natural projection to $X$.

As $z_0$ is also in the center of $\ek \supset \ek_1$, it follows that the complex structure $J$ on $\ep_1$ is $K$-invariant. Using (7), $J$ therefore defines an almost complex structure on $\Fh$. Let $(\, , \,)$ be any $K$-invariant scalar product on $\eg / \ek = \ep_1 \oplus \eg_2 / \ek_2$ and denote by $g$ the corresponding bundle-like metric on $TX$, proposition \ref{t61}.

\begin{prop}
  \label{t63}
$J$ defines a complex structure on $\Fh$. The metric $g$ is K\"ahler-bundle like. In particular $\Fh$ is a K\"ahler--Riemann foliation.
\end{prop}

\begin{proof}
  The leaves of $\Fh$ are isomorphic as almost complex manifolds to \linebreak
$\Gamma_1 \ohne G_1 / K_1$ where $\Gamma_1$ is a discrete subgroup of $G_1$ which acts without fixed points on $G_1 / K_1$. It can depend on the leaf in question. Since $\Gamma_1 \ohne G_1 / K_1$ is a complex manifold, the first assertion follows.

As $(\, , \,)$ is $K_1$-invariant, and $z_0$ is in the center of $\ek_1$ it follows that
\[
g (Jv , w) + g (v, Jw) = 0 \; .
\]
Hence setting:
\[
h (v,w) = g (v,w) + ig (v , Jw) \quad\mbox{for} \; v,w \in T^c \Fh
\]
we get a hermitian metric on $T^c \Fh$. The associated $2$-form along $\Fh$:
\[
\omega_{\Fh} (v,w) = - \halb g (v, Jw)
\]
is closed: Its pullback $\omega_{\Fh_S}$ to $G / K$ is $G$-invariant by construction. Hence the restrictions of $\omega_{\Fh_S}$ to the leaves of $\Fh_S$ i.e. to $G_1 / K_1 \times \{ g_2 K_2 \}$ are left $G_1$-invariant hence closed by a result of E. Cartan. Thus $d_{\Fh_S} \omega_{\Fh_S} = 0$ and therefore $d_{\Fh} \omega_{\Fh} = 0$.
\end{proof}

We will now describe the Hodge structure on the reduced leafwise cohomology of $\Fh$ in terms of the Hodge structure on certain relative Lie algebra cohomologies.

Specializing proposition \ref{t62} and the isomorphisms (10) to the present context we find:
\[
\oH^{\hullet}_{\Fh} (X,\C) = \bar{\bigoplus_{\pi \in {\hat{G}}}} \oH^{\hullet} (\eg_1 , K_1 , H (\pi)^{\infty K_2}) \; .
\]
Now as $G$-representations
\[
H (\pi) = M_{\pi, \Gamma} \otimes V_{\pi}
\]
where $V_{\pi}$ is a fixed representative in the class $\pi$ and where $M_{\pi , \Gamma}$ is the finite dimensional multiplicity space, i.e. 
\[
M_{\pi, \Gamma} = \Hom_G (V_{\pi} , L^2 (\Gamma \ohne G)) \; .
\]
Now since $\hat{G} = \hat{G}_1 \times \hat{G}_2$ we can take $V_{\pi} = V_{\pi_1} \hat{\otimes} V_{\pi_2}$ if $\pi = (\pi_1 , \pi_2)$ and hence
\[
\oH^{\hullet}_{\Fh} (X,\C) = \bar{\bigoplus_{\pi_1 , \pi_2}} M_{\pi , \Gamma} \otimes \oH^{\hullet} (\eg_1 , K_1 , (V_{\pi_1} \hat{\otimes} V^{K_2}_{\pi_2})^{\infty}) \; .
\]
We now assume for simplicity that $G_2$ is also reductive and that $K_2$ is maximal compact in $G_2$. A unitary irreducible representation of a reductive group is admissible by a theorem of Harish-Chandra. In particular $V^{K_2}_{\pi_2}$ and $H^{\hullet} (\eg_1 , K_1 , V^{\infty}_{\pi_1})$ are then (clearly) finite-dimensional, so that we get:
\[
\oH^{\hullet}_{\Fh} (X,\C) = \bar{\bigoplus_{\pi_1 , \pi_2}} M_{\pi , \Gamma} \otimes H^{\hullet} (\eg_1 , K_1 , V^{\infty}_{\pi_1}) \otimes V^{K_2}_{\pi_2} \; .
\]
Here we have also used that $V^{K_2}_{\pi_2}$ consists of smooth vectors. It is clear that in the relative Lie algebra cohomology we can replace $V^{\infty}_{\pi_1}$ by
\[
V_{\pi_1 , 0} = V^{\infty}_{\pi_1} \cap (V_{\pi_1})_{(K_1)}
\]
the intersection being taken with the space $(V_{\pi_1})_{(K)}$ of $K_1$-finite vectors in $V_{\pi_1}$. Then $V_{\pi_1 , 0}$ is a unitary admissible $(\eg , K)$-module in the sense of \cite{BW} II.2.1. 

We now assume the following about our $K$-invariant scalar product $(\,,\,)$ on $\eg / \ek = \ep_1 \oplus \eg_2 / \ek_2$ in order to use results of \cite{BW}: Fix a $G_1$- and $\Theta$-invariant non-degenerate symmetric form $B$ on $\eg_1$ whose restrictions to $\ek_1$ resp. $\ep_1$ are negative resp. positive definite. Choose $(\,,\,)$ such that on $\ep_1$ it agrees with $B$. We define the Casimir element $C  \in \eU (\eg_1)$ as in \cite{BW} II.1.3 using $B$.
It follows from \cite{BW} II.3.1 that
\[
H^{\hullet} (\eg_1 , K_1 , V_{\pi_1 , 0}) = C^{\hullet} (\eg_1 , K_1 , V_{\pi_1,0}) = \Hom_{K_1} (\Lambda^{\hullet} \ep_1 , V_{\pi_1 ,0})
\]
if the central character $\chi_{\pi_1}$ is trivial and if the Casimir element $C$ acts trivially on $V_{\pi_1,0}$.

If one of these conditions is not satisfied, the representation $V_{\pi_1,0}$ has no cohomology:
\[
H^{\hullet} (\eg_1 , K_1 , V_{\pi_1,0}) = 0 \; .
\]
See \cite{VZ} for a classification of representations with non-zero cohomology.

Until now we have therefore seen the following fact:

\begin{theorem}
  \label{t64}
Under the above conditions we have:
\[
\oH^{\hullet}_{\Fh} (X , \C) = \bar{\bigoplus_{\pi_1 , \pi_2}} M_{\pi,\Gamma} \otimes \Hom_{K_1} (\Lambda^{\hullet} \ep_1 , V_{\pi_1 , 0}) \otimes V^{K_2}_{\pi_2} \; .
\]
Here the sum runs over $\pi_1 \in \hat{G}_1 , \pi_2 \in \hat{G}_2$ s.t. $\chi_{\pi_1} = 1 , d \pi_1 (C)= 0$ and $V^{K_2}_{\pi_2} \neq 0$.
\end{theorem}

Note that this result generalizes \cite{BW} VII Cor. 3.4 which treats the case $G_2 = 1$, a result essentially due to Matsushima.

The Hodge decomposition of $\oH^{\hullet}_{\Fh} (X ,\C)$ can now be described using representation theory.

\begin{theorem}
  \label{t65}
Under the isomorphism in theorem \ref{t64} we have
\[
H^{pq} = \bar{\bigoplus_{\pi_1 , \pi_2}} M_{\pi , \Gamma} \otimes \Hom_{K_1} (\Lambda^p \ep^+_1 \otimes \Lambda^q \ep^-_1 , V_{\pi_1,0}) \otimes V^{K_2}_{\pi_2}
\]
with $\pi_1 , \pi_2$ as in \ref{t64}.
\end{theorem}

\begin{proof}
  The assertion is a consequence of the following:

$\bullet$ The Laplacian $\Delta_{\Fh}$ on $\Ah^{\hullet}_{\Fh} (X)$ is identified via the isomorphism
\begin{eqnarray*}
  \Ah^p_{\Fh} (X) & = & C^p (\eg_1 , K_1 , C^{\infty} (\Gamma \ohne G)^{K_2}) \\
& = & \bar{\bigoplus_{\pi_1 , \pi_2}} M_{\pi , \Gamma} \otimes C^p (\eg_1 , K_1 , V_{\pi_1,0}) \otimes V^{K_2}_{\pi_2}
\end{eqnarray*}
with the operator
\[
\bar{\bigoplus_{\pi_1 ,\pi_2}} \id \otimes \Delta_{\pi_1} \otimes V^{K_2}_{\pi_2} \; .
\]
Here $\Delta_{\pi_1}$ is the Laplace operator introduced in \cite{BW} II \S\,2.\\
$\bullet$ The isomorphism:
\begin{eqnarray*}
  \Ah^{pq}_{\Fh} (X) & = & \Hom_{K_1} (\Lambda^p \ep^+_1 \otimes \Lambda^q \ep^-_1 , C^{\infty} (\Gamma \ohne G)^{K_2}) \\
& = & \bigoplus_{\pi_1 , \pi_2} M_{\pi,\Gamma} \otimes \Hom_{K_1} (\Lambda^p \ep^+_1 \otimes \Lambda^q \ep^-_1 , V_{\pi_1,0}) \otimes V^{K_2}_{\pi_2}
\end{eqnarray*}
together with proposition \ref{t44}.
\end{proof}

One can also describe the Lefschetz operator in representation theoretic terms under the isomorphism of theorem \ref{t64}.

\begin{minipage}[t]{7cm}
Mathematisches Institut\\
WWU M\"unster\\
Einsteinstr. 62\\
48149 M\"unster\\
Germany\\
deninge@math.uni-muenster.de
\end{minipage} \hfill
\begin{minipage}[t]{7cm}
Mathematisches Institut\\ 
Universit\"at D\"usseldorf\\
Universit\"atsstr. 1\\
40225 D\"usseldorf\\
Germany\\
singhof@cs.uni-duesseldorf.de
\end{minipage}
\end{document}